\documentclass[11pt]{amsart}
\usepackage[T1]{fontenc}
\usepackage{lmodern}
\usepackage[a4paper,margin=27mm]{geometry}
\usepackage{amsmath,amssymb,amsthm}
\usepackage{microtype}
\usepackage[hidelinks]{hyperref}
\hypersetup{
  pdftitle={Borel (2n+1)-coloring for bounded degree function graphs},
  pdfauthor={Junhao Chen}
}

\newtheorem{theorem}{Theorem}[section]

\newtheorem{question}[theorem]{Question}

\theoremstyle{definition}

\title[Borel $(2n+1)$-coloring for bounded degree function graphs]{Borel $(2n+1)$-coloring for bounded degree function graphs}
\author{Junhao Chen}
\address{School of Mathematical Sciences and LPMC, Nankai University, Tianjin 300071, P.R. China}
\email{jhchen@mail.nankai.edu.cn}

\date{}

\begin{document}

\begin{abstract}
Let $G=(X,E)$ be a Borel graph generated by $n$-many bounded-to-one Borel functions. We show that $G$ admits a Borel proper $(2n+1)$-coloring, which answers a problem of Kechris, Solecki and Todor\v{c}evi\'{c}.
\end{abstract}

\maketitle

\section{Introduction}

Let $X$ be a standard Borel space and let $f_1,\dots,f_n:X\longrightarrow X$ be Borel functions. We say $G=(X,E)$ is generated by $f_1,\ldots,f_n$ if for any two $x,y\in X$,
\begin{center}
$xEy \iff x\ne y$ and for some $i=1,\dots,n$, $f_i(x)=y$ or $f_i(y)=x$.   
\end{center}
When there is no ambiguity, we also say that $G$ is a \textit{function graph}. A proper $k$-coloring is a map $c:X\longrightarrow \{1,\dots,k\}$ such that $c(x)\ne c(y)$ whenever $xEy$. The Borel chromatic number $\chi_B(G)$ is the least cardinality of a standard Borel space $Y$ admitting a Borel proper coloring $c:X\longrightarrow Y$.

The study of Borel coloring of function graphs originates in the work of Kechris, Solecki, and Todor\v{c}evi\'c~\cite{KST99}. For $k\in \mathbb{N}$, we say a function $f:X\longrightarrow X$ is \textit{$\le k$-to-$1$} if for each $x\in X$, $|f^{-1}(\{x\})|\le k$. We say $f$ is \textit{bounded-to-one} if such a $k$ exists.
They proved the following result:

\begin{theorem}[{\cite[Proposition~4.7 and Corollary~5.3]{KST99}}]
Let $X$ be a standard Borel space and $f_1,\dots,f_n:X\longrightarrow X$ be Borel functions, each of which is $\le k$-to-$1$ for some fixed $k$. Then $\chi_B(G)\le \min \{3^n, (k+1)n+1\}$.
\end{theorem}

They also asked the following question, which has been considered in a more general setting in later research.

\begin{question}[{\cite[Problem~4.9]{KST99}; \cite[Problem~5.8]{KM20}}]
Let $G$ be defined as above. Is it true that $\chi_B(G)\le 2n+1$? More generally, for arbitrary Borel functions $f_1,\dots,f_n$, does $\chi_B(G)<\aleph_0$ imply $\chi_B(G)\le 2n+1$?
\end{question}

The following results do not require the functions to be bounded-to-one.
Palamourdas~\cite{Pal12} improved the upper bound to $\frac{1}{2}(n+1)(n+2)$ under the assumption of finite Borel chromatic number. He also obtained the bounds $5$ for two functions and $8$ for three functions. Then Meehan and Palamourdas~\cite[Corollary~2.4]{MP21} improved the general bound to $\frac{1}{2}(n+1)(n+2)-2$ for $n\geq4$. The problem has also been studied in the measurable setting. Pelayo-G\'omez~\cite{PG26} showed that the $\mu$-measurable chromatic number is at most $2n+1$ for a graph generated by $n$-many Borel functions and every Borel probability measure $\mu$. Bernshteyn, Hou, Shinko, and Weilacher~\cite{BHSW26} gave a shorter proof by constructing a random Borel $(2n+1)$-coloring in which each fixed edge is almost surely nonmonochromatic.

In this paper, we answer the question from~\cite[Problem~4.9]{KST99} in the bounded-degree setting.

\begin{theorem}\label{thm:main}
Let $X$ be a standard Borel space and $f_1,\dots,f_n:X\longrightarrow X$ be Borel functions, each of which is $\le k$-to-$1$ for some fixed $k$. Then their function graph admits a  Borel proper $(2n+1)$-coloring.
\end{theorem}

\section{Proof of Theorem~\ref{thm:main}}

\noindent We first fix some notation before the proof.

\noindent Let $D=\{(x,y)\in X^2:x\ne y,\ \exists\,1\leq i\leq n,\ f_i(x)=y\}$.

\noindent For each $x\in X$, let $N^{-}(x)=\bigcup_{1\leq i\leq n}f_i^{-1}(\{x\})\setminus\{x\}$ and $N^{+}(x)=\{f_i(x):1\leq i\leq n\}\setminus\{x\}$.

\noindent For a map $c:X\longrightarrow\{1,\dots,2n+1\}$, we say $x$ is \textit{bad for $c$} if $c(x)=c(y)$ for some $y\in N^{+}(x)$. Otherwise we say $x$ is \textit{good}.
\begin{proof}
Let $\Delta$ be a bound on the degrees of $G$. By~\cite[Proposition~4.6]{KST99}, there is a finite Borel proper coloring $f:X\longrightarrow\{1,\dots,m\}$ for some $m$. We will inductively define a sequence of maps $(c_t)_t$ from $X$ to $\{1,\dots,2n+1\}$ and a sequence of maps $(w_t)_t$ from $D$ to $\mathbb{N}$. The maps $c_t$ are the colorings obtained by recoloring vertices at each stage, and the maps $w_t$ will help us record the informations of these recolorings.

Let $c_0\equiv 1$ and $w_0\equiv 0$. Suppose stage $t-1$ has been defined. At stage $t$, we start with the current state $(c,w)=(c_{t-1},w_{t-1})$ and recolor bad vertices in order of their $f$-colors. For $j=1,\dots,m$, consider each vertex $y$ with $f(y)=j$ that is bad for the current coloring $c$. Let $L_y=\{1,\dots,2n+1\}\setminus\{c(f_i(y)):1\le i\le n,\ f_i(y)\ne y\}$ and note that $|L_y|\geq n+1$. Then for each $a\in L_y$, let $s_a=\sum_{x\in N^{-}(y),\,c(x)=a}w(x,y)$.

We choose the least $l\in L_y$ such that $s_l=\min\{s_a:a\in L_y\}$ and set $c(y)=l$. For each $x\in N^{-}(y)$, we add $n$ to $w(x,y)$ when $c(x)=l$, and replace it by $\max\{w(x,y)-1,0\}$ otherwise. After all $m$ classes have been processed, we set $(c_t,w_t)=(c,w)$.

We now estimate the number of recolorings up to stage $t$.  Consider a single recoloring of $y$ up to this stage, with chosen color $l$. Let $\Phi_y=\sum_{x\in N^{-}(y)}w(x,y)^2$ and $W=\sum_{x\in N^{-}(y)}w(x,y)$, where $c$ and $w$ denote the coloring and weights before this recoloring.
By the choice of $l$, we have:
\[
\begin{aligned}
(n+1)s_l \leq |L_y|s_l\leq\sum_{a\in L_y}s_a
 =\sum_{\substack{x\in N^{-}(y)\\c(x)\in L_y}}w(x,y)\leq\sum_{x\in N^{-}(y)}w(x,y)=W.
\end{aligned}
\] 
Put $A=\{x\in N^{-}(y):c(x)=l\}$ and $B=N^{-}(y)\setminus A$. Then $\sum_{x\in A}w(x,y)=s_l$ and $\sum_{x\in B}w(x,y)=W-s_l$. For each $x\in A$, the weight $w(x,y)$ increases by $n$. So the contribution of part $A$ is:
\[
\sum_{x\in A}\bigl((w(x,y)+n)^2-w(x,y)^2\bigr)=2ns_l+n^2|A|.
\]
For each $x\in B$, the weight $w(x,y)$ is replaced by $\max\{w(x,y)-1,0\}$. If $w(x,y)\geq1$, the corresponding change in its square is $-2w(x,y)+1$; if $w(x,y)=0$, the change is zero and is still at most $-2w(x,y)+1$. Thus the contribution of part $B$ is at most:
\[
\sum_{x\in B}\bigl(-2w(x,y)+1\bigr)=-2(W-s_l)+|B|.
\]
Let $\Phi_{y}^{'}$ denote the value of $\Phi_y$ after this recoloring. Adding the two contributions and using $(n+1)s_l\leq W$, we obtain
\[
\begin{aligned}
\Phi_{y}^{'}-\Phi_y &\leq 2ns_l+n^2|A|-2(W-s_l)+|B|=2(n+1)s_l-2W+n^2|A|+|B|\\
&\leq n^2(|A|+|B|)=n^2|N^{-}(y)|\leq n^2\Delta.
\end{aligned}
\]

From now on, $w$ and $\Phi_y$ denote their values at the end of stage $t$. For each $z\in X$, let $\varphi(z)$ denote the number of times $z$ has been recolored from the start of the construction up to stage $t$. For $1\leq j\leq\varphi(y)$, let $\ell_j$ be the new color assigned to $y$ at its $j$-th recoloring. For each $x\in N^{-}(y)$, let $a_{x,j}$ be the color of $x$ immediately before that recoloring, and define
\[
h(x,y)=\left|\left\{j\in\{1,\dots,\varphi(y)\}:\ell_j=a_{x,j}\right\}\right|.
\]
In other words, $h(x,y)$ counts the number of recolorings of $y$ up to stage $t$ in which the new color of $y$ equals the color of $x$ at that moment.
Since $x$ and $y$ have different $f$-colors, $x$ is not recolored at the same time as $y$. The weight $w(x,y)$ increases by $n$ on these $h(x,y)$ occasions and decreases by at most one on the remaining $\varphi(y)-h(x,y)$ occasions. Since the weight is initially zero, we have
\[
w(x,y)\geq nh(x,y)-(\varphi(y)-h(x,y))=(n+1)h(x,y)-\varphi(y).
\]
On the other hand, $w(x,y)^2$ is one of the nonnegative terms in $\Phi_y$. Since the weights are initially zero and the weights $w(x,y)$ change only when $y$ is recolored, summing the bound $\Phi_y'-\Phi_y\leq n^2\Delta$ over these $\varphi(y)$ recolorings gives $\Phi_y\leq n^2\Delta\varphi(y)$. Thus
\[
w(x,y)\leq\sqrt{\Phi_y}\leq n\sqrt{\Delta\varphi(y)}.
\]
Combining these two bounds, we obtain
\[
h(x,y)\leq\frac{\varphi(y)+n\sqrt{\Delta\varphi(y)}}{n+1}.
\]

Let $M=\max_{x\in X}\varphi(x)\leq t$. If $M>0$, choose $x$ with $\varphi(x)=M$. After each recoloring, $x$ is good. Before its next recoloring, some $y\in N^{+}(x)$ must  be recolored with the color of $x$, which remains unchanged between these two recolorings. So at least $M-1$ such recolorings must have occurred. Hence we have:
\[
M-1\leq\sum_{y\in N^{+}(x)}h(x,y).
\]

For each $y\in N^{+}(x)$, applying the preceding estimate and using $\varphi(y)\leq M$, we obtain
\[
h(x,y)\leq\frac{\varphi(y)+n\sqrt{\Delta\varphi(y)}}{n+1}
\leq\frac{M+n\sqrt{\Delta M}}{n+1}.
\]
Hence
\[
M\leq n+1+n^2\sqrt{M\Delta}.
\]
Since $2n^2\sqrt{M\Delta}\leq M+n^4\Delta$, we obtain $M\leq n^4\Delta+2n+2$, independently of $t$. Every vertex is therefore recolored only finitely often. Let $c=\lim_{t\to\infty}c_t$. If a vertex $x$ were bad for $c$, then once its color and those of its out-neighbors $N^+(x)$ had stabilized, it would have to be recolored again, a contradiction. Thus $c$ is a Borel proper $(2n+1)$-coloring.
\end{proof}

At the end of this article we do some remarks. Our proof use a weighting function $w$ to control the number of recoloring. Forte Shinko suggests that we can define a more explicit recoloring rule and proceed a more simple combinatorial counting argument. These two ways led to the same result--- the color of each vertex will be stable. However, we also note that the final color of each vertex may not be decided by information of a finite radius ball. In fact, we do not know  in which stage the recoloring happens. In other words, our solution may not be \textit{continuous}, with the setting that  the underlying space $X$ is a zero-dimensional Polish space and $\mathbb{N}$ equipped with discrete topology. This is also known as a study of \textit{Continuous Combinatorics}. For more details on continuous colorings and \textbf{\textsf{LOCAL}} algorithms, see, e.g., \cite[Section~2.B.2]{Ber23a}. And see
\cite{GJ15,GJKS25,Ber23b,BCGGRV22,GR23,GWW25} for recent progress on Continuous Combinatorics.

\section*{Acknowledgement}
We would like to thank Forte Shinko for suggesting the more elementary combinatorial argument. And we thank Anton Bernshteyn and Felix Weilacher for pointing out the relationship between our algorithm and continuous combinatorics in zero-dimensional Polish space. We also thank Edward Hou for many helpful conversations throughout this work.

\section*{AI Disclosure}

We used ChatGPT 6 Pro to explore a deterministic adaptation of the randomized recoloring argument of Bernshteyn, Hou, Shinko, and Weilacher~\cite{BHSW26}. These discussions led to an AI-generated draft of the proof, and the author wrote the whole article based on this draft. We take full responsibility for the contents and correctness of the article.

\end{document}